\UseRawInputEncoding
\documentclass[a4paper, 11pt]{article}

\usepackage{authblk}
\usepackage{setspace}
\usepackage{comment}
\usepackage{amssymb}
\usepackage{amsthm}
\usepackage{amscd}
\usepackage{graphicx}
\usepackage{titlesec}
\usepackage{indentfirst}
\usepackage{cite}
\titleformat{\section}{\bfseries}{\thesection .}{0.5em}{}
\titleformat{\subsection}{\itshape}{\thesubsection .}{0.5em}{}
\usepackage[top=25.4mm, bottom=25.4mm, left=27.7mm, right=27.8mm]{geometry}

\usepackage[tbtags]{amsmath}
\usepackage{booktabs}
\allowdisplaybreaks  
\usepackage{mathrsfs}
\usepackage{caption}
\begin{document}

\newtheorem{theorem}{Theorem}[section]  
\newtheorem{example}{Example}[section]            
\newtheorem{algorithm}{Algorithm}[section]

\newtheorem{axiom}{Axiom}[section]
\newtheorem{property}{Property}[section]

\newtheorem{proposition}[theorem]{Proposition}%
\newtheorem{corollary}[theorem]{Corollary}%
\newtheorem{definition}[theorem]{Definition}%
\newtheorem{lemma}{Lemma}[section]
\newtheorem{cor}{Corollary}[section]
\newtheorem{remark}{Remark}[section]
\newtheorem{condition}{Condition}[section]
\newtheorem{conclusion}{Conclusion}[section]
\newtheorem{assumption}{Assumption}[section]

\title{\Large Trotter-Kato approximations of semilinear stochastic evolution equations in Hilbert spaces}
\author{
{\normalsize Xia Zhang\thanks{
\textit{Email addresses}: zhangxia@tiangong.edu.cn (Xia Zhang)}}\\
{\normalsize Lingfei Dai\thanks{
\textit{Email addresses}: dlf970226@163.com (Lingfei Dai)}}\\
{\normalsize Ming Liu\thanks{
\textit{Email addresses}: liuming@tiangong.edu.cn (Ming Liu)}}\\

\textit{\small School of Mathematical Sciences, Tiangong University, Tianjin 300387, P.R.China}}




\date{}
\maketitle

\renewcommand{\baselinestretch}{1.2}
\large\normalsize
\noindent \rule[0.5pt]{14.5cm}{0.6pt}\\
\noindent
\textbf{Abstract}\\
Motivated by the work of T.E. Govindan in [5,8,9], this paper is concerned with a more general semilinear stochastic evolution equation.
The difference between the equations considered in this paper and the previous one is that it makes some changes to the nonlinear function in random integral, which also depends on the probability distribution of stochastic process at that time.
First, this paper considers the existence and uniqueness of mild solutions for such equations.
Furthermore, Trotter-Kato approximation system is introduced for the mild solutions, and the weak convergence of induced probability measures and zeroth-order approximations are obtained.
Then we consider the classical limit theorem about the parameter dependence of this kind of equations.
Finally, an example of stochastic partial differential equation is given to illustrate our results.

\vspace{0.3cm}

\noindent \rule[0.5pt]{14.5cm}{0.6pt}\\
\noindent
\textbf{Mathematics Subject Classification (2021 MSC):} 60H15
\section{Introduction}\label{sec1}
Problems related to mild solutions of stochastic evolution equations have been studied by many scholars, for instance, Ahmed and Ding \cite{ANUDX95}, Da Prato and Zabczyk \cite{DPGZJ92}, Govindan and Ahmed \cite{GTEANU13}, Govindan \cite{GTE16}, Ichikawa \cite{IA82} and McKean \cite{MHP66} etc. Such problems are often used to simulate the stochastic process which are found in the study of dynamic system in many fields such as science, engineering and finance.

In view of such issues, the previous scholars have done a lot of excellent work, and the relative results are widely used in real life, which are very worthy of promotion. Therefore, this is also the significance of this paper. In this paper, the nonlinear function term corresponding to the stochastic integral term is further generalized on the basis of predecessors.

The aim of this paper is to study the stochastic process $\{x(t),\ t\in[0,\infty)\}$ described by a semilinear stochastic evolution equation in a real separable Hilbert space $H_1$ of the form as follows:
\begin{align}
\label{system(1.1)}
  dx(t)&=[Ax(t)+f(x(t),\mu(t))]dt+g(x(t),\mu(t))d\omega(t),\quad t\in(0,\infty)\\
\label{system(1.2)}
  x(0)&=x_0,
 \end{align}
 where $A$ is the infinitesimal generator of a $C_0$-semigroup $\{S(t),\ t\in[0,\infty)\}$ of bounded linear operators on $H_1;$ $f$ is an $H_1$-valued function defined on $H_1\times V_{\varphi^2}(H_1),$ where $V_{\varphi^2}(H_1)$ is a subset of probability measures of $H_1$; $\mu(t)$ is the probability distribution of $x(t)$ for all $t\in(0,\infty);$ $g$ is a $\mathcal{L}(H_2,H_1)$-valued function defined on $H_1\times V_{\varphi^2}(H_1)$; $\omega(t)$ is a $H_2$-valued Wiener process; and the initial value $x_0$ is a $\mathcal{F}_0$-measurable $H_1$-valued random variable.

In Da Prato and Zabczyk \cite{DPGZJ92}, it has been systematically explained that when the nonlinear drift term $f(t,x)$ and the constant additive diffusion term $g(t,x)$ only depend on the state of solution process $x(t)$ at that time. The specific equation form is shown as follows:
\begin{align}
\label{system(1.3)}
  dx(t)&=[Ax(t)+f(t,x)]dt+g(t,x)d\omega(t),\quad t\in(0,\infty)\\
\label{system(1.4)}
  x(0)&=x_0,
 \end{align}
 in Ichikawa \cite{IA82}, Yosida approximation system for this kind of equations is introduced and its specific application is given. In Govindan \cite{GTE15}, Trotter-Kato approximation system is introduced and its specific application is obtained. In Govindan \cite{GTE20}, the weak convergence of the probability measures corresponding to the mild solutions of such equations in Trotter Kato approximation system is proved.

In Ahmed and Ding \cite{ANUDX95}, it has studied the case when the specific diffusion term $g(x)=\sqrt{Q}$. That is, the equation becomes the following form:
\begin{align}
\label{system(1.5)}
  dx(t)&=[Ax(t)+f(x(t),\mu(t))]dt+\sqrt{Q}d\omega(t),\quad t\in(0,\infty)\\
\label{system(1.6)}
  x(0)&=x_0,
 \end{align}
where $Q$ is a positive, symmetric and bounded operator on $H_1$. Ahmed and Ding proved the existence and uniqueness of mild solutions of Equation (\ref{system(1.5)}), and introduced the Yosida approximation system. Then in Govindan \cite{GTE06}, this paper introduced and studied the Trotter-Kato approximation system for the System (\ref{system(1.5)})-(\ref{system(1.6)}).

In 2015, Govindan and Ahmed \cite{GTEANU15} considered the case when the nonlinear drift term $f$ depends on the state and probability distribution $\mu(t)$ of the solution process $x(t)$ at that time. The form of the equation is as follows:
\begin{align}
\label{system(1.7)}
  dx(t)&=[Ax(t)+f(x(t),\mu(t))]dt+g(x(t))d\omega(t),\quad t\in(0,\infty)\\
\label{system(1.8)}
  x(0)&=x_0,
 \end{align}
Govindan and Ahmed studied the existence and uniqueness of mild solutions, and introduced Yosida approximation system for the System (\ref{system(1.7)})-(\ref{system(1.8)}). Furthermore, the weak convergence of the induced probability measures corresponding to mild solutions in the system is obtained, and the exponential stability of mild solutions as a practical application in the system is given. In 2018, the Trotter-Kato approximation system was introduced and studied in Govindan \cite{GTE18}, then the weak convergence of the induced probability measures of mild solutions and specific application in this system were obtained.

It can be seen from the above that the constant additive diffusion term $g$ depends on the state and probability distribution $\mu(t)$ of the solution process $x(t)$ at that time, which has not been considered. Hence, this is also the main motivation of this paper. The specific form of equation is shown in System (\ref{system(1.1)})-(\ref{system(1.2)}). Firstly, we show that the existence and uniqueness of mild solutions for such equations; then we introduce Trotter-Kato approximation system, consider the zeroth-order approximation and other related problems (see Kannan and Barucha-Reid \cite{KDBRAT85} and Govindan \cite{GTE94}), and give the approximation error estimates; finally, we give the classical limit theorem about the parameter dependence of equation as a specific application (see Gikhman and Skorokhod \cite{GIISAV72}).

 The main difficulty of this paper is that the change of the constant additive diffusion term $g$ will lead to many previously proved conclusions that need to be re-proven. For example, the existence and uniqueness of mild solutions of such equations, and various conclusions of equations in the Trotter-Kato approximation system.

 The remainder of this paper is organized as follows: in section 2, we will give the preliminary knowledge; in section 3, we will give a concrete explanation of the existence and uniqueness of mild solutions for the System (\ref{system(1.1)})-(\ref{system(1.2)}); in Section 4, we introduce the Trotter-Kato approximation system for this kind of equations; in section 5, we will study the parameter dependence of the equations in Trotter-Kato approximation system.

\section{Preliminaries}\label{sec2}
In the sequel, the terms we use mainly derived from Govindan \cite{GTE18}.
Let $H_1$ and $H_2$ be real separable Hilbert spaces and $\mathcal{L}(H_2,H_1)$ the space of bounded linear operators from $H_2$ to $H_1$. We use $|\cdot|$ and $(\cdot,\cdot)$ to represent the norms and scalar products of $H_1$ and $H_2$. We write $\mathcal{L}(H_1)$ for $\mathcal{L}(H_1,H_1).$ Let $(\Omega,\mathcal{F},P)$ be a complete probability space. A map $x:\Omega\rightarrow H_1$ is a random variable if it is strongly measurable. Let $x:\Omega\rightarrow H_1$ be a square integrable random variable, namely, $x\in L_2(\Omega,\mathcal{F},P;H_1)\equiv L_2(\Omega;H_1).$ The covariance operator of the random variable $x$ is $Cov[x]=E[(x-Ex)\circ(x-Ex)]$, where $E$ denotes the expectation and for any $h_1,h_2\in\mathcal{L}(H_1),$ we have $h_1\circ h_2\in\mathcal{L}(H_1)$ is defined by $(h_1\circ h_2)h_3=h_1(h_2,h_3),$ $h_3\in H_1.$ Then $Cov[x]$ is a nonnegative selfadjoint trace class operator and satisfies tr$Cov[x]=E[x-Ex]^2$, where tr denotes the trace, and trace class operator is also called nuclear operator. The joint covariance of any pair $\{x,h\}\in L_2(\Omega;H_1)$ is defined by $Cov[x,h]=E[(x-Ex)\circ(h-Eh)].$

Let $T\in(0,\infty).$ The stochastic process $\{x(t),\ t\in[0,T]\}$ is a family of random variables with values of $H_1.$ Let $\{\mathcal{F}_t,\ t\in[0,T]\}$ be a family of increasing sub-$\sigma$-algebras of the $\sigma$-algebra $\mathcal{F}$. If $x(t)$ is $\mathcal{F}_t$-measurable for all $t\in[0,T]$, then the stochastic process $\{x(t),\ t\in[0,T]\}$ is said to be $\mathcal{F}_t$-adapted.

The stochastic process $\{\omega(t),\ t\in[0,T]\}$ in $H_2$ is a $Q$-Wiener process if $\omega(t)\in L_2(\Omega;H_2)$ and $E\omega(t)=0$ for all $t\in[0,T];$ $Cov[\omega(t)-\omega(u)]=(t-u)Q,\ Q\in L_1^+(H_2)$ is a nonnegative trace class operator; $\omega(t)$ has continuous sample paths and independent increments. The operator $Q$ is called the incremental covariance operator of Wiener process $\omega(t).$ Therefore, $\omega$ has the expression $\omega(t)=\sum_{n=1}^\infty\eta_n(t)e_n,$ where $\{e_n,n\in\mathbb{N^+}\}$ is an orthonormal set of eigenvectors of $Q,$ $\eta_n(t),\ n\in\mathbb{N^+}$, are mutual independent real-valued Wiener processes with incremental covariance $\kappa_n>0$, $Qe_n=\kappa_ne_n$ and $trQ=\sum_{n=1}^\infty\kappa_n.$

Let $\{S(t),\ t\in[0,T]\}$ be a $C_0$-semigroup of bounded linear operators on $H_1$ and $A$ is the infinitesimal generator of $\{S(t),\ t\in[0,T]\}.$ There exist constants $\alpha\geqslant0$ and $M\geqslant1$ such that
$$\|S(t)\|\leqslant Me^{\alpha t}\quad t\in[0,T],$$
where $\|\cdot\|$ denotes the operator norm in $\mathcal{L}(H_1).$ We use the notation $A\in G(M,\alpha)$ to denote this concept.

Let $\mathcal{B}(H_1)$ denote the Borel $\sigma$-algebra on $H_1$, and $M(H_1)$ denote the space composed of probability measures with the usual topology of weak convergence on $\mathcal{B}(H_1)$. $C(H_1)$ denotes the space composed of all continuous functions mapped from $H_1$ to $H_1.$ We use $(\mu,\psi)$ to represent $\int_{H_1}\psi(x)\mu(dx)$ when this integral makes sense. We let $\varphi(x)\equiv1+|x|,$ $x\in H_1$, then we define the Banach space (see Govindan and Ahmed \cite{GTEANU15}):
$$C_\psi(H_1)=\bigg\{\psi\in C(H_1):\|\psi\|_{C_\psi(H_1)}\equiv\mathop{{\rm sup}}_{x\in H_1}\frac{\psi(x)}{\varphi^2(x)}+\mathop{{\rm sup}}_{x\neq h\in H_1}\frac{|\psi(x)-\psi(h)}{|x-h|}<\infty\bigg\}.$$
For $p\geqslant1,$ let $V_{\varphi^p}^u(H_1)$ be the Banach space composed of signed measures $v$ on $H_1$ satisfying $\|\mu\|_{\varphi^p}\equiv\int_{H_1}\varphi^p|v|(dx)<\infty,$ where $|v|=v^++v^-,$ $v^+$ is the positive part of $v$ and $v^-$ is the negative part of $v.$ Let $V_{\varphi^2}(H_1)=V_{\varphi^2}^u(H_1)\cap M(H_1)$ be the set of probability measures with second moments on $\mathcal{B}(H_1).$ We place a topology on $V_{\varphi^2}(H_1)$ induced by the metric as follows:
$$\rho(\mu,\nu)=\sup\bigg\{(\psi,\mu-\nu):\|\psi\|_\rho=\mathop{{\rm sup}}_{x\in H_1}\frac{\psi(x)}{\varphi^2(x)}+\mathop{{\rm sup}}_{x\neq h\in H_1}\frac{|\psi(x)-\psi(h)}{|x-h|}\leqslant1\bigg\}.$$
Hence, $(V_{\varphi^2}(H_1),\rho)$ is a complete metric space. We use $C([0,T],(V_{\varphi^2}(H_1),\rho))$ to represent a complete metric space composed of all continuous functions from $[0,T]$ to $(V_{\varphi^2}(H_1),\rho).$ The metric for this space is $D_t(\mu,\nu)=\mathop{{\rm sup}}_{t\in[0,T]}\rho(\mu(t),\nu(t))$, for $\mu,\nu\in C([0,T],(V_{\varphi^2}(H_1),\rho)).$

Let $C([0,T],L_2(\Omega;H_1))$ be the Banach space consisting of all continuous functions from $[0,T]$ to $L_2(\Omega;H_1)$ satisfying the condition $\mathop{{\rm sup}}_{t\in[0,T]}E|x(t)|^2<\infty.$ Let $\Theta$ be the closed subspace of $C([0,T],L_2(\Omega;H_1))$ composed of $\mathcal{F}_t$-adapted processes $x=\{x(t),\ t\in[0,T]\}.$ Thus $\Theta$ is also a Banach space and its norm is $$\|x\|_\Theta=\sqrt{\mathop{{\rm sup}}_{t\in[0,T]}E|x(t)|^2}<\infty.$$

All stochastic processes mentioned in this paper are supposed to be based on the complete filtered probability space $(\Omega,\mathcal{F},\{\mathcal{F}_t,\ t\in[0,T]\},P)$.

Furthermore, we also need to introduce the following definitions.
\begin{definition}
\label{def2.1}
An $H_1$-valued stochastic process $x=\{x(t),\ t\in[0,T]\}$ with probability law $P(x(t))=\mu(t)$ is said to be a strong solution of Equation (\ref{system(1.1)}) if
\begin{spacing}{1.5}
\end{spacing}
\noindent
$(a)$ for each $t\in[0,T],$ $x(t)$ is $\mathcal{F}_t$-adapted;
\begin{spacing}{1.5}
\end{spacing}
\noindent
$(b)$ $t\rightarrow x(t)$ is continuous $P-a.s.$;
\begin{spacing}{1.5}
\end{spacing}
\noindent
$(c)$ for each $t\in[0,T],$ $x(t)\in D(A)$ (domain of $A$), $P-a.s.$ and $\int_0^T|Ax(t)|^2dt<\infty,$ $P-a.s.$,
\par\setlength\parindent{1.5em}
 for any $T\in(0,\infty)$ and
 $$x(t)=x_0+\int_0^t[Ax(s)+f(x(s),\mu(s))]ds+\int_0^tg(x(s),\mu(s))d\omega(s),\quad t\in[0,T],\ P-a.s..$$
\end{definition}
It clear that the conditions of this definition are rather strong, thus we introduce the following definition.
\begin{definition}
\label{def2.2}
An $H_1$-valued stochastic process $x=\{x(t),\ t\in[0,T]\}$ with probability law $P(x(t))=\mu(t)$ is called a mild solution of System (\ref{system(1.1)})-(\ref{system(1.2)}), or only Equation (\ref{system(1.1)}) if
\begin{spacing}{1.5}
\end{spacing}
\noindent
$(a)$ $x$ is jointly measurable and $\mathcal{F}_t$-adapted and its restriction to the interval $[0,T]$ satisfies
\par\setlength\parindent{1.5em}
$\int_0^t|x(t)|^2dt<\infty,$ $P-a.s.$;
\begin{spacing}{1.5}
\end{spacing}
\noindent
$(b)$ $x(t)$ satisfies the following integral equation
\begin{equation*}
\begin{split}
x(t)&=S(t)x_0+\int_0^tS(t-s)f(x(s),\mu(s))ds+\int_0^tS(t-s)g(x(s),\mu(s))d\omega(s),\quad P-a.s.
\end{split}
\end{equation*}
\par\setlength\parindent{1.5em}
for all $t\in[0,T].$
\end{definition}
The stochastic integral term in the last equality is defined in the sense of it\^{o} (see Da Prato and Zabczyk \cite{DPGZJ92}).
\begin{proposition}
\label{por(2.3)}
\rm{({Ichikawa \cite{IA82} Proposition 1.9})} \it Let $M(H_1,H_2)$ be the space of stochastic process $G(\cdot,\cdot):[0,T]\times\Omega\rightarrow\mathcal{L}(H_1,H_2)$ which are strongly measurable. Let $G\in M(H_1,H_2)$ with $\int_0^TE|G(t)|^pd(t)<\infty$ some integer $p\geqslant2,$ and let $y(t)=\int_0^tG(r)d\omega(r).$ Then
\begin{equation*}
\begin{split}
E|y(t)|^p&\leqslant[\frac{1}{2}p(p-1)]^{\frac{p}{2}}\left[\int_0^t\left[E(trG(r)WG^\ast(r))^{\frac{p}{2}}\right]^{\frac{2}{p}}dr\right]^{\frac{p}{2}}\\
&\leqslant[\frac{1}{2}p(p-1)]^{\frac{p}{2}}(trW)^{\frac{p}{2}}t^{\frac{p}{2}-1}\int_0^tE|G(r)|^pdr,\quad t\in[0,T].
\end{split}
\end{equation*}
\end{proposition}

\section{Existence and Uniqueness of a Mild Solution}\label{sec3}

In this section, we state and prove the existence and uniqueness of mild solutions of System (\ref{system(1.1)})-(\ref{system(1.2)}). This is the first main result of this paper. Here we mainly consider the case that the nonlinear term $g$ also depends on the probability distribution of mild solutions $x$ at that time.

First of all, we make the following assumption:

\begin{spacing}{2.0}
\begin{large}
\noindent
\textbf{Assumption(A1)}
\end{large}
\end{spacing}

(i) Let $\{S(t),\ t\in[0,T]\}$ be a $C_0$-semigroup of bounded linear operators on $H_1$ and $A$ is
\par\setlength\parindent{1.5em}
the infinitesimal generator of $\{S(t),\ t\in[0,T]\};$
\begin{spacing}{1.5}
\end{spacing}

(ii) For $p \geqslant 2$, $f : H_1\times (V_{\varphi^{2}}(H_1),\rho)\rightarrow H_1$ and $g : H_1\times (V_{\varphi^{2}}(H_1),\rho)\rightarrow \mathcal{L}(H_2,H_1)$ satisfy
\par\setlength\parindent{1.5em}
the following Lipschitz and linear growth conditions :
$$|f(x,\mu)-f(y,\nu)| \leqslant K_{1}(|x-y|+\rho(\mu,\nu)),$$
$$|g(x,\mu)-g(y,\nu)| \leqslant K_{2}(|x-y|+\rho(\mu,\nu)),$$
$$|f(x,\mu)|^{p} \leqslant K_{3}(1+|x|^{p}+\|\mu\|^{^{p}}_{\varphi}),$$
$$|g(x,\mu)|^{p} \leqslant K_{4}(1+|x|^{p}+\|\mu\|^{^{p}}_{\varphi}),$$

\begin{spacing}{1.5}
\end{spacing}
for all $x, y \in H_1$ and $\mu, \nu \in V_{\varphi^{2}}(H_1)$, where $K_{i}, i=1,2,3,4$ are positive numbers and
\par\setlength\parindent{1.5em}
independent of time $t$.

\begin{theorem}
\label{the3.1}
Assume that the Assumption (A1) hold. Then for each initial value, which is $\mathcal{F}_{0}$-measurable $H_1$-valued random variable $x_{0} \in L_{2}(\Omega;H_1)$,
\begin{spacing}{1.5}
\end{spacing}

(a) The System (\ref{system(1.1)})-(\ref{system(1.2)}) has a unique mild solution $x=\{x(t),\ t\in[0,T]\}$ in $\Theta$ and

\qquad its corresponding probability distribution $\mu=\{\mu(t)=P(x(t)),\ t\in[0,T]\}$ in

\qquad$C([0,T], (V_{\varphi^{2}}(H_1),\rho))$.

\begin{spacing}{1.5}
\end{spacing}
(b) For any $p \geqslant 1$ and $\mathcal{F}_{0}$-measurable initial value $x_{0} \in L_{2p}(\Omega;H_1)$, we have the

\qquad following inequality:
$$\mathop{{\rm sup}}_{t\in[0,T]}E|x(t)|^{2p}\leqslant J(1+E|x_0|^{2p}),$$
\qquad\  where $J$ is a positive number.
\end{theorem}

\begin{proof}
\begin{spacing}{1.5}
\end{spacing}
First, we use the Banach contraction mapping principle and Picard's iterations to show that for every fixed measure-valued function $\mu \in C([0,T], (V_{\varphi^{2}}(H_1),\rho))$, System (\ref{system(1.1)})-(\ref{system(1.2)}) has a unique mild solution $x_{\mu}=\{x_{\mu}(t):\ t\in [0,T]\}$ in $\Theta$. For $\mu \in C([0,T], (V_{\varphi^{2}}(H_1),\rho))$, let's define the operator $\Lambda_\mu$ on $\Theta$ by the following integral equation:
$$(\Lambda_{\mu}x)(t)=S(t)x_0+\int_0^t S(t-s)f(x(s),\mu(s))ds+\int_0^t S(t-s)g(x(s),\mu(s))d \omega(s),$$
for $t \in [0,T]$. According to Assumption (A1), we have
\begin{equation*}
\begin{split}
&E|(\Lambda_\mu x)(t)|^2\\
&\leqslant3\bigg\{E|S(t)x_0|^2+E\left|\int_0^tS(t-s)f(x(s),\mu(s))ds\right|^2+E\left|\int_0^tS(t-s)g(x(s),\mu(s))d\omega(s)\right|^2\bigg\}\\
&\leqslant3\bigg\{M^2e^{2\alpha T}E|x_0|^2+TM^2e^{2\alpha T}K_3\int_0^t(1+E|x(s)|^2+\|\mu(s)\|_\varphi^2)ds\\
&\quad+(trQ)M^2e^{2\alpha T}K_4\int_0^t(1+E|x(s)|^2+\|\mu(s)\|_\varphi^2)ds\bigg\}\\
\end{split}
\end{equation*}
\begin{equation*}
\begin{split}
&\leqslant3M^2e^{2\alpha T}E|x_0|^2+[3T^2M^2e^{2\alpha T}K_3+3(trQ)M^2e^{2\alpha T}K_4T]\left(1+\|x\|_\Theta^2+\mathop{{\rm sup}}_{t\in[0,T]}\|\mu(t)\|_\varphi^2\right)\\
&<\infty.
 \end{split}
\end{equation*}
Hence, $\mathop{{\rm sup}}_{t\in[0,T]}E|(\Lambda_\mu x)(t)|^2<\infty$. It is clear that $(\Lambda_\mu x)(t)$ is $\mathcal{F}_t$-measurable when $x(t)$ is $\mathcal{F}_t$-measurable. Therefore, $\Lambda_\mu$ maps $\Theta$ into $\Theta$, that is $\|\Lambda_\mu x\|_{\Theta}<\infty.$

\begin{spacing}{1.5}
\end{spacing}
Next, for any $x,y \in \Theta$ with $x(0)=y(0)$ (namely, the initial values of $x$ and $y$ are the same), from the Assumption (A1), we get
\begin{align}
  E|(\Lambda_\mu x)(t)-(\Lambda_\mu y)(t)|^2&\leqslant2\bigg\{E\left|\int_0^tS(t-s)[f(x(s),\mu(s))-f(y(s),\mu(s))]ds\right|^2\nonumber\\
  &\quad+E\left|\int_0^tS(t-s)[g(x(s),\mu(s))-g(y(s),\mu(s))]d\omega(s)\right|^2\bigg\}\nonumber\\
  &\leqslant2TM^2e^{2\alpha T}\int_0^tE\left|f(x(s),\mu(s))-f(y(s),\mu(s))\right|^2ds\nonumber\\
  &\quad+2(trQ)M^2e^{2\alpha T}\int_0^tE\left|g(x(s),\mu(s))-g(y(s),\mu(s))\right|^2ds\nonumber\\
  &\leqslant4TM^2e^{2\alpha T}{K_1}^2\int_0^t(E|x(s)-y(s)|^2+\rho^2(\mu(s),\mu(s)))ds\nonumber\\
  &\quad+4(trQ)M^2e^{2\alpha T}{K_2}^2\int_0^t(E|x(s)-y(s)|^2+\rho^2(\mu(s),\mu(s)))ds\nonumber\\
  &=4M^2e^{2\alpha T}[T{K_1}^2+(trQ){K_2}^2]\int_0^tE|x(s)-y(s)|^2ds.\label{eq(3.1)}
 \end{align}

 For any $n\in\mathbb{N^+}$, from the Picard's iterations, when Equation (\ref{eq(3.1)}) is iterated $n$ times, we get
\begin{equation*}
\begin{split}
&\mathop{{\rm sup}}_{t\in[0,T]}E|(\Lambda_\mu^n x)(t)-(\Lambda_\mu^n y)(t)|^2\\
&\leqslant4M^2e^{2\alpha T}[T{K_1}^2+(trQ){K_2}^2]\int_0^tE|\Lambda_\mu^{n-1}x(s)-\Lambda_\mu^{n-1}y(s)|^2ds.\\
\end{split}
\end{equation*}
Then we obtain
 $$\|\Lambda_\mu^n x-\Lambda_\mu^n y\|_{\Theta}^2\leqslant\frac{\left\{4M^2e^{2\alpha T}[T{K_1}^2+(trQ){K_2}^2]\right\}^nT^n}{n!}\|x-y\|_\Theta^2.$$

When $n$ gets big enough, we have $\left\{4M^2e^{2\alpha T}[T{K_1}^2+(trQ){K_2}^2]\right\}^nT^n/{n!}<1$, according to Banach contraction mapping principle, $\Lambda_\mu^n$ is a contraction map on $\Theta$. Hence, from Picard's iteration, we know that $\Lambda_\mu$ has a unique fixed point in $\Theta.$ Then for a given $\mu\in C([0,T],(V_{\varphi^2}(H_1),\rho))$, the System (\ref{system(1.1)})-(\ref{system(1.2)}) has a unique mild solution $x_\mu$.

\begin{spacing}{1.5}
\end{spacing}
 Next, let $P(x_\mu)=\left\{P(x_\mu(t)):t\in[0,T]\right\}$ denote the probability distribution of mild solutions $x_\mu.$ Now we prove that $P(x_\mu)\in C([0,T],(V_{\varphi^2}(H_1),\rho)).$ Since we already know $x_\mu\in \Theta,\ P(x_\mu(t))\in V_{\varphi^2}(H_1)$ for all $t\in [0,T],$ it follows that we only need to prove $t\rightarrow P(x_\mu(t))$ is continuous. According to the semigroup property of $S(t)$, it is clear that for $0\leqslant s\leqslant t\leqslant T,$ we obtain
\begin{align*}
  x_\mu(t)-x_\mu(s)&=[S(t)-S(s)]x_0+\int_0^t[S(t-u)-S(s-u)]f(x(u),\mu(u))du\nonumber\\
  &\quad+\int_0^t[S(t-u)-S(s-u)]g(x(u),\mu(u))d\omega(u)\nonumber\\
  &=S(s)[S(t-s)-I]x_0+\int_s^tS(t-u)f(x(u),\mu(u))du\nonumber\\
  &\quad+\int_0^s[S(t-s)-I]S(s-u)f(x(u),\mu(u))du\\
  &\quad+\int_s^tS(t-u)g(x(u),\mu(u))d\omega(u)\nonumber\\
  &\quad+\int_0^s[S(t-s)-I]S(s-u)g(x(u),\mu(u))d\omega(u).\nonumber
\end{align*}
Then taking the expectations on both sides of the above expression we arrive at the following inequality:
\begin{equation}
\label{eq(3.2)}
\begin{split}
  E\left| x_\mu(t)-x_\mu(s)\right|^2&\leqslant5E\left|S(s)[S(t-s)-I]x_0\right|^2+5E\left|\int_s^tS(t-u)f(x(u),\mu(u))du\right|^2\\
  &\quad+5E\left|\int_0^s[S(t-s)-I]S(s-u)f(x(u),\mu(u))du\right|^2\\
  &\quad+5E\left|\int_s^tS(t-u)g(x(u),\mu(u))d\omega(u)\right|^2\\
  &\quad+5E\left|\int_0^s[S(t-s)-I]S(s-u)g(x(u),\mu(u))d\omega(u)\right|^2.
\end{split}
\end{equation}

Now let's estimate each term on the RHS of (\ref{eq(3.2)}): since $S(t)$ is a $C_0$-semigroup, it follows that for any $x\in H_1,$ $[S(t-s)-I]x$ converges to 0 as $t\rightarrow s.$ We consider the first and the third term. We already know $x_0\in L_2(\Omega;H_1),\ f(x(t),\mu(t))\in L_2(\Omega;H_1)$ for $t\in [0,T]$ a.e.. According to Assumption (A1), we have
$$E|S(s)[S(t-s)-I]x_0|^2\leqslant M^4e^{4\alpha T}E|x_0|^2<\infty,$$
and
\begin{flalign*}
\begin{split}
&E\left|\int_0^s[S(t-s)-I]S(s-u)f(x(u),\mu(u))du\right|^2\\
&\leqslant TM^4e^{4\alpha T}\int_0^sE\left|f(x(u),\mu(u))\right|^2du\\
&\leqslant T^2M^4e^{4\alpha T}K_3(1+\|x\|_{\Theta}^2+\mathop{{\rm sup}}_{t\in[0,T]}\|\mu(t)\|_\varphi^2)\\
&<\infty.
\end{split}
\end{flalign*}
From the strong continuity of $S(t)$ and the Lebesgue's dominated convergence theorem, we have
$$\mathop{{\rm lim}}_{t\rightarrow s}E|S(s)[S(t-s)-I]x_0|^2=0,$$
and
$$\mathop{{\rm lim}}_{t\rightarrow s}E\left|\int_0^s[S(t-s)-I]S(s-u)f(x(u),\mu(u))du\right|^2=0.$$
Next we consider the fifth term in Equation (\ref{eq(3.2)}). For given $x\in\Theta,$ it follows from the Assumption (A1) that $g(x(t))\in L_2([0,T];\mathcal{L}(H_2,H_1)).$ Hence, for the stochastic integrals term, also from the strong continuity of semigroup we get
\begin{equation*}
\begin{split}
&E\left|\int_0^s[S(t-s)-I]S(s-u)g(x(u),\mu(u))d\omega(u)\right|^2\\
&\leqslant (trQ)M^4e^{4\alpha T}\int_0^sE\left|g(x(u),\mu(u))\right|^2du\\
&\leqslant (trQ)M^4e^{4\alpha T}K_4T(1+\|x\|_{\Theta}^2+\mathop{{\rm sup}}_{t\in[0,T]}\|\mu(t)\|_\varphi^2)\\
&<\infty.\\
 \end{split}
\end{equation*}
By the Lebesgue's dominated convergence theorem, we have
$$\mathop{{\rm lim}}_{t\rightarrow s}E\left|\int_0^s[S(t-s)-I]S(s-u)g(x(u),\mu(u))d\omega(u)\right|^2=0.$$
We consider the second and the fourth term in Equation (\ref{eq(3.2)}),
\begin{equation*}
\begin{split}
 &\mathop{{\rm lim}}_{t\rightarrow s}E\left|\int_s^tS(t-u)f(x(u),\mu(u))du\right|^2\\
 &\leqslant (t-s)M^2e^{2\alpha T}\int_s^tE\left|f(x(u),\mu(u))\right|^2du\\
 &\leqslant (t-s)^2M^2e^{2\alpha T}K_3(1+\|x\|_{\Theta}^2+\mathop{{\rm sup}}_{t\in[0,T]}\|\mu(t)\|_\varphi^2)\\
 &=0,
\end{split}
\end{equation*}
and
\begin{equation*}
\begin{split}
 &\mathop{{\rm lim}}_{t\rightarrow s}E\left|\int_s^tS(t-u)g(x(u),\mu(u))d\omega(u)\right|^2\\
 &\leqslant(trQ)M^2e^{2\alpha T}\int_s^tE\left|g(x(u),\mu(u))\right|^2du\\
 &\leqslant (trQ)M^2e^{2\alpha T}K_4(t-s)(1+\|x\|_{\Theta}^2+\mathop{{\rm sup}}_{t\in[0,T]}\|\mu(t)\|_\varphi^2)\\
 &=0.
\end{split}
\end{equation*}
Finally, we have
$$\mathop{{\rm lim}}_{t\rightarrow s}E|x_\mu(t)-x_\mu(s)|^2=0.$$
For any $\psi$ in $C_\psi(H_1),$ from the definition of the metric $\rho$ and Jensen's inequality, we get
\begin{equation*}
\begin{split}
|(\psi,P(x_\mu)(t)-P(x_\mu)(s))|&\leqslant E|\psi(x_\mu(t))-\psi(x_\mu(s))|\\
&\leqslant\|\psi\|_\rho E|x_\mu(t)-x_\mu(s)|\\
&\leqslant E|x_\mu(t)-x_\mu(s)|\rightarrow0\quad as\ t\rightarrow s,
\end{split}
\end{equation*}
thus we obtain
$$\rho(P(x_\mu)(t),P(x_\mu)(s))\rightarrow0\quad as\ t\rightarrow s.$$
From this we know that $P(x_\mu)(t)\in C([0,T],(V_{\varphi^2}(H_1),\rho)).$
\begin{spacing}{1.5}
\end{spacing}
Now we prove that the probability law $\mu$ is also unique for every mild solution $x_\mu$. First, we need to define an operator $\Delta$ on $C([0,T],(V_{\varphi^2}(H_1),\rho))$ as follows:
$$\Delta:\mu\rightarrow P(x_\mu).$$
Assume that $x$ is a mild solution of System (\ref{system(1.1)})-(\ref{system(1.2)}). Then it is obvious that the probability law $P(x)=\mu$ is a fixed point of $\Delta.$ On the contrary, if $\mu$ is a fixed point of $\Delta,$ then the probability law $\mu$ corresponding to the mild solution $x_\mu$ that defined by the parameterization of $\mu$ in System (\ref{system(1.1)})-(\ref{system(1.2)}) also belongs to $C([0,T],(V_{\varphi^2}(H_1),\rho))$. Therefore, we only need to prove that the operator $\Delta$ has a unique fixed point in $C([0,T],(V_{\varphi^2}(H_1),\rho)).$

For this, let $\mu$ and $\nu$ be any two elements of $C([0,T],(V_{\varphi^2}(H_1),\rho))$ and let $x_\mu$ and $x_\nu$ be the mild solutions of the parameterization of their respective probability laws with the same initial values. These two mild solutions satisfy the integral equations from the Definition \ref{def2.2} (b), respectively. Then consider the difference between the two integral equations:
\begin{equation*}
\begin{split}
x_\mu(t)-x_\nu(t)&=\int_0^tS(t-s)[f(x_\mu(s),\mu(s))-f(x_\nu(s),\nu(s))]ds\\
&\quad+\int_0^tS(t-s)[g(x_\mu(s),\mu(s))-g(x_\nu(s),\nu(s))]d\omega(s)
\end{split}
\end{equation*}
Next, taking the expectations on both sides of the expression, by the Assumption (A1), we obtain
\begin{equation*}
\begin{split}
 &E|x_\mu(t)-x_\nu(t)|\\
 &\leqslant2\bigg\{E\left|\int_0^tS(t-s)[f(x_\mu(s),\mu(s))-f(x_\nu(s),\nu(s))]ds\right|^2\\
 &\quad+E\left|\int_0^tS(t-s)[g(x_\mu(s),\mu(s))-g(x_\nu(s),\nu(s))]d\omega(s)\right|^2\bigg\}\\
 &\leqslant2TM^2e^{2\alpha T}\int_0^tE|f(x_\mu(s),\mu(s))-f(x_\nu(s),\nu(s))|^2ds\\
 &\quad+2(trQ)M^2e^{2\alpha T}\int_0^tE|g(x_\mu(s),\mu(s))-g(x_\nu(s),\nu(s))|^2ds\\
 &\leqslant4M^2e^{2\alpha T}[TK_1^2+(trQ)K_2^2]E\int_0^t[|x_\mu(s)-x_\nu(s)|^2+\rho^2(\mu(s),\nu(s))]ds.\\
\end{split}
\end{equation*}
Hence, we have
\begin{equation*}
\begin{split}
&\mathop{{\rm sup}}_{t\in[0,T]}E|x_\mu(t)-x_\nu(t)|^2\\
&\leqslant4TM^2e^{2\alpha T}[T{K_1}^2+(trQ)K_2^2]\bigg\{\mathop{{\rm sup}}_{t\in[0,T]}E|x_\mu(t)-x_\nu(t)|^2+D_t(\mu,\nu)\bigg\}.\\
\end{split}
\end{equation*}
Given $T\in(0,\infty)$, let $T$ be small enough to satisfy the inequality as follows:
$$4TM^2e^{2\alpha T}[T{K_1}^2+(trQ){K_2}^2]<\frac{1}{3}.$$
Considering the above two inequalities, we get
\begin{equation}
\label{eq(3.3)}
\begin{split}
\mathop{{\rm sup}}_{t\in[0,T]}E|x_\mu(t)-x_\nu(t)|^2&<\frac{1}{2}D_t^2(\mu,\nu).
\end{split}
\end{equation}
Then from the Equation (\ref{eq(3.3)}) and the inequality $\rho^2(P(x_\mu(t)),P(x_\nu(t)))\leqslant E|x_\mu(t)-x_\nu(t)|^2,$ we have
$$D_t^2(\Delta(\mu),\Delta(\nu))<\frac{1}{2}D_t^2(\mu,\nu).$$
Therefore, let small $T=T',$ $\Delta$ is a contraction map in $C([0,T'],(V_{\varphi^2}(H_1),\rho))$, then $\Delta$ has a unique fixed point. Hence, the System (\ref{system(1.1)})-(\ref{system(1.2)}) has a unique solution $x$ with its corresponding probability distribution $\mu\equiv P(H_1)$ on $[0,T']$. Since $x\in\Theta$, it follows that the solution can be extended to any finite interval $[0,T]$ by extended to $[T',2T'],\ [2T',3T'],\ [3T',4T'],$ etc. we have completed the proof of part (a).

For any $p\geqslant 1,$ by the Assumption (A1), we get
\begin{align}
\label{eq(3.4)}
  E|x(t)|^{2p}&\leqslant 3^{2p}\bigg\{E|S(t)x_0|^{2p}+E\left|\int_0^tS(t-s)f(x(s),\mu(s))ds\right|^{2p}\nonumber\\
  &\quad+E\left|\int_0^tS(t-s)g(x(s),\mu(s))d\omega(s)\right|^{2p}\bigg\}\nonumber\\
  &\leqslant3^{2p}\bigg\{M^{2p}E|x_0|^{2p}+T^{2p-1}M^{2p}K_3\int_0^t(1+E|x(s)|^{2p}+\|\mu(s)\|_\varphi^{2p})ds\nonumber\\
  &\quad+(trW)^pM^{2p}K_4\int_0^t(1+E|x(s)|^{2p}+\|\mu(s)\|_\varphi^{2p})ds\bigg\}\nonumber\\
  &\leqslant j_1+j_2E|x_0|^{2p}+j_3\int_0^tE|x(s)^{2p}ds,\qquad t\in [0,T],
\end{align}
where $j_i,\ i=1,2,3$ are positive numbers. The Gronwall's inequality is applied to Equation (\ref{eq(3.4)}), we get
\begin{equation*}
\begin{split}
E|x(t)|^{2p}&\leqslant(j_1+j_2E|x_0|^{2p})\exp(\int_0^tj_3ds)\\
&\leqslant j_1\exp(j_3T)+j_2\exp(j_3T)E|x_0|^{2p}，
\end{split}
\end{equation*}
let $J=\max\{j_1\exp(j_3T);\ j_2\exp(j_3T)\}$, then we have
$$\mathop{{\rm sup}}_{t\in[0,T]}E|x(t)|^{2p}\leqslant J(1+E|x_0|^{2p}),$$
where $J$ depends on the parameters that appear in Equation (\ref{eq(3.4)}). Now the proof is completed.
\end{proof}

Next let us deal with the problem of continuous dependence on the initial conditions of mild solutions. This problem can be solved by the conclusion of Theorem \ref{the3.1}.

Let $x(t)$ and $y(t)$ for $t\in[0,T]$ be the two mild solutions of Equation (\ref{system(1.1)}), and their corresponding initial values are $x_0$ and $y_0$, respectively. For all $t\in[0,T]$, let $\mu(t)=P(x(t))$ and $\nu(t)=P(y(t))$ denote the corresponding probability distributions.
\begin{corollary}
\label{cor3.2}
Suppose that the conditions of Theorem \ref{the3.1} are hold, then for each $T\in[0,\infty)$, there exists a constant $C$ such that
$$\mathop{{\rm sup}}_{0\leqslant s\leqslant t}E|x(s)-y(s)|^2\leqslant CE|x_0-y_0|^2 \qquad t\in[0,T],$$
where $C=3M^2e^{2\alpha T}\exp\{12TM^2e^{2\alpha T}[TK_1^2+(trQ)K_2^2]\}$.
\end{corollary}
\begin{proof}
First, we differentiate the two mild solutions
\begin{equation*}
\begin{split}
x(s)-y(s)&=S(s)x_0-S(s)y_0+\int_0^sS(s-u)[f(x(u),\mu(u))-f(y(u),\nu(u))]du\\
&\quad+\int_0^sS(s-u)[g(x(u),\mu(u))-g(y(u),\nu(u))]d\omega(u).
\end{split}
\end{equation*}
let us square this equation and then take the expectations on both sides, we arrive at the following inequality:
\begin{equation*}
\begin{split}
&E|x(s)-y(s)|^2\\
&\leqslant3\bigg\{E\left|S(s)(x_0-y_0)\right|^2+E\left|S(s-u)[f(x(u),\mu(u))-f(y(u),\nu(u))]du\right|^2\\
&\quad+E\left|\int_0^sS(s-u)[g(x(u),\mu(u))-g(y(u),\nu(u))]d\omega(u)\right|^2\bigg\}\\
&\leqslant3M^2e^{2\alpha T}E|x_0-y_0|^2+6TM^2e^{2\alpha T}K_1^2\int_0^s(E|x(u)-y(u)|^2+\rho^2(\mu(u),\nu(u)))du\\
&\quad+6(trQ)M^2e^{2\alpha T}K_2^2\int_0^s(E|x(u)-y(u)|^2+\rho^2(\mu(u),\nu(u)))du\\
&\leqslant3M^2e^{2\alpha T}E|x_0-y_0|^2+12M^2e^{2\alpha T}[TK_1^2+(trQ)K_2^2]\int_0^sE|x(u)-y(u)|^2du,
\end{split}
\end{equation*}
where $\rho^2(\mu(u),\nu(u))\leqslant E|x(u)-y(u)|^2,$ then from the Gronwall's inequality we get the following inequality:
$$\mathop{{\rm sup}}_{0\leqslant s\leqslant t}E|x(s)-y(s)|^2\leqslant CE|x_0-y_0|^2 \qquad t\in[0,T].$$
Thus the corollary is proved.
\end{proof}

\section{Trotter-Kato Approximations}\label{sec4}
In this section, we will set up the Trotter-Kato Approximations system for the System (\ref{system(1.1)})-(\ref{system(1.2)}) and estimate the error in the approximation.

First let's study the following family of stochastic evolution equations:
\begin{align}
\label{system(4.1)}
  dx_n(t)&=[A_nx_n(t)+f(x_n(t),\mu_n(t))]dt+g(x_n(t),\mu_n(t))d\omega(t),\quad t>0\\
\label{system(4.2)}
  x_n(0)&=x_0,
 \end{align}
 where $A_n,\ n\in\mathbb{N^+},$ is the infinitesimal generator of a $C_0$-semigroup $\{S_n(t),\ t\geqslant0\}$ of bounded linear operators on $H_1.$

 From the Theorem \ref{the3.1} (a), it follows that for every $n\in\mathbb{N^+}$, the System (\ref{system(4.1)})-(\ref{system(4.2)}) has a unique mild solution $x_n\in C([0,T], L_2(\Omega;H_1)),$ then $x_n(t)$ satisfies the following stochastic integral equation:
 \begin{equation*}
\begin{split}
x_n(t)&=S_n(t)x_0+\int_0^tS_n(t-s)f(x_n(s),\mu_n(s))ds\\
&\quad+\int_0^tS_n(t-s)g(x_n(s),\mu_n(s))d\omega(s),\quad t\in[0,T],\quad P-a.s.
\end{split}
\end{equation*}

In order to introduce the Trotter-Kato Approximation system, we first need to make the following assumptions:
\begin{spacing}{2.0}
\begin{large}
\noindent
\textbf{Assumption(A2)}
\end{large}
\end{spacing}
(i)$A_n\in G\left(M,\alpha\right)$ for every $n\in\mathbb{N^+}$;
\begin{spacing}{1.5}
\end{spacing}
(ii)As $n\rightarrow\infty,$ $A_nx\rightarrow Ax$ for every $x\in D,$ where $D$ is a dense subset of $X$;
\begin{spacing}{1.5}
\end{spacing}
(iii)There exist a  $\xi$ with $Re\ \xi>\alpha$ for which $(\xi I-A)D$ is dense in $X$, then the closure

\ \ \ \ \ $\overline{A}$ of $A$ is in $G\left(M,\alpha\right)$.
\begin{spacing}{1.5}
\end{spacing}
\begin{proposition}
\label{the4.1}
$\rm(Pazy\ [Theorem\ 4.5, p.88])$\it Suppose that the Assumption (A2) hold. If $S_n(t)$ and $S(t)$ are the $C_0$-semigroups by $A_n$ and $\overline{A}$, respectively, then
\begin{equation}
\begin{split}
\label{eq(4.3)}
\mathop{{\rm lim}}_{n\rightarrow\infty}S_n(t)x&=S(t)x, \qquad x\in X,
\end{split}
\end{equation}
for all $t\geqslant0,$ and the limit in Equation (\ref{system(4.1)}) is uniform in $t$ for $t$ in bounded intervals.
\end{proposition}
\begin{theorem}
\label{the4.2}
Assume that the Assumption (A1) and (A2) are satisfied. Let $x(t)$ and $x_n(t)$ be the mild solutions of Equations (\ref{system(1.1)}) and (\ref{system(4.1)}), respectively. Then for every $T\in(0,\infty),$ we have
$$\mathop{{\rm lim}}_{n\rightarrow\infty}\mathop{{\rm sup}}_{0\leqslant t\leqslant T}E|x_n(t)-x(t)|^2=0.$$
\end{theorem}
\begin{proof}
First we take the difference between the two mild solutions

\begin{equation*}
\begin{split}
&x_n(t)-x(t)\\
&=[S_n(t)-S(t)]x_0+\int_0^t[S_n(t-s)f(x_n(s),\mu_n(s))-S(t-s)f(x(s),\mu(s))]ds\\
&\quad+\int_0^t[S_n(t-s)g(x_n(s),\mu_n(s))-S(t-s)g(x(s),\mu(s))]d\omega(s),\quad P-a.s.,
\end{split}
\end{equation*}
for $t\in[0,T]$. Second let us square this equation and then take the expectations on both sides, then we have

\begin{equation}
\begin{split}
\label{eq(4.4)}
&E|x_n(t)-x(t)|^2\\
&\leqslant5\bigg\{|E[S_n(t)-S(t)]x_0|^2+E\left|\int_0^tS_n(t-s)[f(x_n(s),\mu_n(s))-f(x(s),\mu(s)))]ds\right|^2\\
&\quad+E\left|\int_0^t[S_n(t-s)-S(t-s)]f(x(s),\mu(s))ds\right|^2\\
&\quad+E\left|\int_0^tS_n(t-s)[g(x_n(s),\mu_n(s))-g(x(s),\mu(s)))]d\omega(s)\right|^2\\
&\quad+E\left|\int_0^t[S_n(t-s)-S(t-s)]g(x(s),\mu(s))d\omega(s)\right|^2\bigg\},\quad P-a.s.
\end{split}
\end{equation}
for $t\in [0,T].$

Now let us consider each of the right hand side terms in Equation (\ref{eq(4.4)}). First let's consider the second and fourth term, from the Assumption (A1) and (A2), we have
\begin{align}
\label{eq(4.5)}
&\mathop{{\rm sup}}_{0\leqslant s\leqslant t}E\left|\int_0^sS_n(s-u)[f(x_n(u),\mu_n(u))-f(x(u),\mu(u)))]du\right|^2\nonumber\\
&\leqslant TM^2e^{2\alpha T}K_1^2\int_0^t[E|x_n(s)-x(s)|^2+\rho^2(\mu_n(s),\mu(s))]ds\nonumber\\
&\leqslant 2TM^2e^{2\alpha T}K_1^2\int_0^tE|x_n(s)-x(s)|^2ds.
\end{align}
Next, let us study the fourth term. By Proposition \ref{por(2.3)}, we obtain
\begin{align}
\label{eq(4.6)}
&\mathop{{\rm sup}}_{0\leqslant s\leqslant t}E\left|\int_0^sS_n(s-u)[g(x_n(u),\mu_n(u))-g(x(u),\mu(u)))]d\omega(u)\right|^2\nonumber\\
&\leqslant (trQ)\int_0^t\|S_n(s-u)\|^2E|g(x_n(u),\mu_n(u))-g(x(u),\mu(u))|^2du\nonumber\\
&\leqslant 2(trQ)M^2e^{2\alpha T}K_2^2\int_0^tE|x_n(s)-x(s)|^2ds.
\end{align}
According to Equations (\ref{eq(4.5)}) and (\ref{eq(4.6)}), it is clear that Equation (\ref{eq(4.4)}) can be reduced that
\begin{equation*}
\begin{split}
&\mathop{{\rm sup}}_{0\leqslant s\leqslant t}E|x_n(s)-x(s)|^2\\
&\leqslant\sigma(n,t)+10M^2e^{2\alpha T}[TK_1^2+(trQ)K_2^2]\int_0^tE|x_n(s)-x(s)|^2ds,
\end{split}
\end{equation*}
it is clear that
\begin{equation}
\begin{split}
\label{eq(4.7)}
\sigma(n,t)&=5\mathop{{\rm sup}}_{0\leqslant s\leqslant t}\bigg\{|E[S_n(t)-S(t)]x_0|^2+E\left|\int_0^s[S_n(s-u)-S(s-u)]f(x(u),\mu(u))du\right|^2\\
&\quad+E\left|\int_0^s[S_n(s-u)-S(s-u)]g(x(u),\mu(u))d\omega(u)\right|^2\bigg\}.
\end{split}
\end{equation}
Obviously, $\sigma(n,t)\geqslant0,$ then from the Gronwall's inequality we obtain
$$\mathop{{\rm sup}}_{0\leqslant s\leqslant t}E|x_n(s)-x(s)|^2\leqslant\sigma(n,t)\exp\{10tM^2e^{2\alpha T}[TK_1^2+(trQ)K_2^2]\}, \quad t\in[0,T].$$
Therefore, we only need to prove $\sigma(n,t)\rightarrow0$ as $n\rightarrow\infty.$ We consider the first term of Equation (\ref{eq(4.7)}). Since $\overline{A}\in G(M,\alpha)$ and $A_n\in G\left(M,\alpha\right)$ for every $n\in\mathbb{N^+}$, and
$$E|[S_n(t)-S(t)]x_0|^2\leqslant2M^2e^{2\alpha T}E|x_0|^2<\infty,$$
uniformly in $n$ and $t\in[0,T]$. Hence, from the Theorem \ref{the4.1} and the Lebesgue's dominated convergence theorem, we obtain
\begin{align}
\label{eq(4.8)}
\mathop{{\rm lim}}_{n\rightarrow\infty}\mathop{{\rm sup}}_{0\leqslant t\leqslant T}E|[S_n(t)-S(t)]x_0|^2=0
\end{align}
for all $t\in[0,T],\ x_0\in H_1,$ and the limit in Equation (\ref{eq(4.8)}) is uniform in $t$. Hence, by the Equation (\ref{eq(4.8)}) we have the first term of Equation (\ref{eq(4.7)}) tends to zero as $n\rightarrow\infty$. Next we study the second term of Equation (\ref{eq(4.7)}), by Assumption (A1), (A2) and the Theorem \ref{the3.1}\ (b), we obtain
\begin{equation*}
\begin{split}
&\mathop{{\rm sup}}_{0\leqslant s\leqslant t}E\left|\int_0^s[S_n(s-u)-S(s-u)]f(x(u),\mu(u))du\right|^2\\
&\leqslant2TM^2e^{2\alpha T}K_3\int_0^t(1+E|x(u)|^2+\|\mu(u)\|_\varphi^2)du\\
&\leqslant2T^2M^2e^{2\alpha T}K_3[1+J(1+E|x_0|^2)+\mathop{{\rm sup}}_{0\leqslant t\leqslant T}\|\mu(t)\|_\varphi^2]\\
&<\infty.
\end{split}
\end{equation*}
Thus by the Lebesgue's dominated convergence theorem, we know that the second term of Equation (\ref{eq(4.7)}) also tends to zero as $n\rightarrow\infty.$
Then we consider the third term of Equation (\ref{eq(4.7)}),
\begin{equation*}
\begin{split}
&\mathop{{\rm sup}}_{0\leqslant s\leqslant t}E\left|\int_0^s[S_n(s-u)-S(s-u)]g(x(u),\mu(u))d\omega(u)\right|^2\\
&\leqslant2(trW)M^2e^{2\alpha T}K_4\int_0^t(1+E|x(u)|^2+\|\mu(u)\|_\varphi^2)du\\
&\leqslant2T(trW)M^2e^{2\alpha T}K_4[1+J(1+E|x_0|^2)+\mathop{{\rm sup}}_{0\leqslant t\leqslant T}\|\mu(t)\|_\varphi^2]\\
&<\infty.
\end{split}
\end{equation*}
Finally, the third term also tends to zero as $n\rightarrow\infty$. Hence, $\sigma(n,t)\rightarrow0$ as $n\rightarrow\infty.$ Then we get the following conclusion
$$\mathop{{\rm lim}}_{n\rightarrow\infty}\mathop{{\rm sup}}_{0\leqslant s\leqslant t}E|x_n(s)-x(s)|^2=0.$$
The proof of this theorem is completed.
\end{proof}
\begin{corollary}
The family of probability distributions $\{\mu_n\}_{n=1}^\infty$ corresponding to the mild solutions $\{x_n\}_{n=1}^\infty$ of Equation (\ref{system(4.1)}) converges to the probability distribution $\mu$ of mild solutions $x$ of Equation (\ref{system(1.1)}) in $C([0,T],V_{\varphi^2}(H_1),\rho)$, as $n\rightarrow\infty.$
\end{corollary}
\begin{proof}
From the foregoing fact and Theorem \ref{the4.2}, we get
$$D_t(\mu_n,\mu)=\mathop{{\rm sup}}_{0\leqslant t\leqslant T}\rho(\mu_n(t),\mu(t))\leqslant\mathop{{\rm sup}}_{0\leqslant t\leqslant T}(E|x_n(t)-x(t)|^2)^\frac{1}{2}\rightarrow0,\quad as\ n\rightarrow\infty.$$
The proof is completed.
\end{proof}

Next we consider the zeroth-order approximations. In other words, stochastic evolution equation can be approximated by a deterministic evolution equation.

Let us study the following stochastic evolution equation:
\begin{align}
\label{system(4.9)}
  dx_\varepsilon(t)&=[A_\varepsilon x_\varepsilon(t)+f(x_\varepsilon(t),\mu_\varepsilon(t))]dt+\varepsilon g(x_\varepsilon(t),\mu_\varepsilon(t))d\omega(t),\quad t\in[0,T],\\
\label{system(4.10)}
  x_\varepsilon(0)&=x_0\in D(A_\varepsilon),
 \end{align}
 where $A_\varepsilon(\varepsilon>0)$ is the infinitesimal generator of a $C_0$-semigroup $\{S_\varepsilon(t),\ t\geqslant0\}$ of bounded linear operators on $H_1$. the deterministic evolution equation as follows:
 \begin{align}
\label{system(4.11)}
  dx'(t)&=Ax'(t)+f(x'(t),\mu(t)dt,\quad t\in[0,T],\\
\label{system(4.12)}
  x'(0)&=x_0\in D(A).
 \end{align}
 The mild solutions of Equations (\ref{system(4.9)}) and (\ref{system(4.11)}) satisfy the following integral equations respectively,
 \begin{equation}
\begin{split}
\label{eq(4.13)}
x_\varepsilon(t)&=S_\varepsilon(t)x_0+\int_0^tS_\varepsilon(t-s)f(x_\varepsilon(s),\mu_\varepsilon(s))ds\\
&\quad+\varepsilon\int_0^tS_\varepsilon(t-s)g(x_\varepsilon(s),\mu_\varepsilon(s))d\omega(s),\quad P-a.s.,
\end{split}
\end{equation}
for $t\in[0,T]$, and
 \begin{equation}
\begin{split}
\label{eq(4.14)}
x'(t)&=S(t)x_0+\int_0^tS(t-s)f(x'(s),\mu(s))ds\quad t\in[0,T].
\end{split}
\end{equation}
For any $\varepsilon>0$ , from Theorem \ref{the3.1} we know that the Equation (\ref{system(4.9)}) has a unique mild solution $x_\varepsilon$ in $C([0,T],L_2(\Omega;H_1)),$ the specific form of the solution is given by Equation (\ref{eq(4.13)}); and Equation (\ref{system(4.11)}) also has a unique mild solution, its form is given by Equation (\ref{eq(4.14)}), which is a special case of $g\equiv0.$

For further research, we refer to Kannan and Bharucha-Reid \cite{KDBRAT85} and make the following assumption.
\begin{spacing}{2.0}
\begin{large}
\noindent
\textbf{Assumption(A3)}
\end{large}
\end{spacing}
Let $A, A_\varepsilon\in G(M,\alpha)(\varepsilon>0)$ with $D(A_\varepsilon)=D(A);\ \mathop{{\rm lim}}_{\varepsilon\rightarrow0}S_\varepsilon(t)=S(t)$ uniformly in $t\in[0,T]$ for each $T\in(0,\infty)$.

Now we show that how to estimate the error of the approximation.
\begin{theorem}
Assume that the Assumption (A1) and (A3) hold. Let $x_\varepsilon(t)$ and $x'(t)$ satisfy the corresponding integral equations. Then
$$E|x_\varepsilon(t)-x'(t)|^2\leqslant\tau(\varepsilon)\eta(t),$$
where $\tau(\varepsilon)$ is a positive monotone decreasing function to zero as $\varepsilon\downarrow0$; $\eta(t)$ is a positive exponentially increasing function.
\end{theorem}
\begin{proof}
 First we take the difference between the two mild solutions,
\begin{equation}
\begin{split}
\label{eq(4.15)}
&x_\varepsilon(t)-x'(t)\\
&=[S_\varepsilon(t)-S(t)]x_0+\int_0^t[S_\varepsilon(t-s)f(x_\varepsilon(s),\mu_\varepsilon(s))-S(t-s)f(x'(s),\mu'(s))ds\\
&\quad+\varepsilon\int_0^tS_\varepsilon(t-s)g(x_\varepsilon(s),\mu_\varepsilon(s))d\omega(s)\\
&=[S_\varepsilon(t)-S(t)]x_0+\int_0^t[S_\varepsilon(t-s)[f(x_\varepsilon(s),\mu_\varepsilon(s))-f(x'(s),\mu'(s))]ds\\
&\quad+\int_0^t[S_\varepsilon(t-s)-S(t-s)]f(x'(s),\mu(s))ds+\varepsilon\int_0^tS_\varepsilon(t-s)g(x'(s),\mu'(s))d\omega(s)\\
&\quad+\varepsilon\int_0^tS_\varepsilon(t-s)[g(x_\varepsilon(s),\mu_\varepsilon(s))-g(x'(s),\mu'(s))d\omega(s)\quad P-a.s.
\end{split}
\end{equation}
for $t\in[0,T].$

Now we start estimating each of terms on the right hand side of Equation (\ref{eq(4.15)}).

Since $\mathop{{\rm lim}}_{\varepsilon\rightarrow0}S_\varepsilon(t)=S(t)$ uniformly in $t\in[0,T]$, it follows that there exists an $\varepsilon_1>0$ and a constant $l_1>0$ such that for any $\delta_1(\varepsilon)>0,$ $\varepsilon$ is dependent only on $\delta_1(\varepsilon)$ and independent of $t$, when $\varepsilon_1>\varepsilon$, we have $E\left|[S_\varepsilon(t)-S(t)]x_0\right|^2\leqslant l_1\delta_1(\varepsilon)$, for all $t\in[0,T].$

Now let us study the second term of Equation (\ref{eq(4.15)}), we have
\begin{equation*}
\begin{split}
&E\left|\int_0^tS_\varepsilon(t-s)[f(x_\varepsilon(s),\mu_\varepsilon(s))-f(x'(s),\mu'(s))]ds\right|^2\\
&\leqslant TM^2e^{2\alpha T}\int_0^tE\left|f(x_\varepsilon(s),\mu_\varepsilon(s))-f(x'(s),\mu'(s))\right|^2ds\\
&\leqslant4TM^2e^{2\alpha T}K_1^2\int_0^tE|x_\varepsilon(s)-x'(s)|^2ds.
\end{split}
\end{equation*}
Consider the third term as well, we obtain
\begin{equation*}
\begin{split}
&E\left|\int_0^t[S_\varepsilon(t-s)-S(t-s)]f(x'(s),\mu'(s))ds\right|^2\\
&\leqslant2TM^2e^{2\alpha T}\int_0^tE\left|f(x'(s),\mu'(s))\right|^2ds\\
&\leqslant2TM^2e^{2\alpha T}K_3\int_0^t(1+E|x'(s)|^2+\|\mu'(s)\|_\varphi^2)ds\\
&\leqslant2T^2M^2e^{2\alpha T}K_3[1+J(1+E|x_0|^2)+\mathop{{\rm sup}}_{t\in[0,T]}\|\mu'(t)\|_\varphi^2]\\
&<\infty.
\end{split}
\end{equation*}
Thus, from the Lebesgue's dominated convergence theorem, we get
$$E\left|\int_0^t[S_\varepsilon(t-s)-S(t-s)]f(x'(s),\mu'(s))ds\right|^2\rightarrow0\quad as\ \ \varepsilon\downarrow0.$$
As in the case of the first term, for any $\delta_2(\varepsilon)>0,$ there exists an $\varepsilon_2>\varepsilon\downarrow0$ and a constant $l_2>0$ such that
$$E\left|\int_0^t[S_\varepsilon(t-s)-S(t-s)]f(x'(s),\mu'(s))ds\right|^2<l_2\delta_2(\varepsilon).$$
Finally, let us consider two stochastic integral terms,
\begin{equation*}
\begin{split}
&E\left|\varepsilon\int_0^tS_\varepsilon(t-s)[g(x_\varepsilon(s),\mu_\varepsilon(s))-g(x'(s),\mu'(s))]d\omega(s)\right|^2\\
&\leqslant\varepsilon^2M^2e^{2\alpha T}(trQ)\int_0^tE|g(x_\varepsilon(s),\mu_\varepsilon(s))-g(x'(s),\mu'(s))|^2ds\\
&\leqslant2\varepsilon^2M^2e^{2\alpha T}(trQ)K_2^2\int_0^t[E|x_\varepsilon(s)-x'(s)|^2+\rho^2(\mu_\varepsilon(s),\mu'(s))]ds\\
&\leqslant4\varepsilon^2M^2e^{2\alpha T}(trQ)K_2^2\int_0^tE|x_\varepsilon(s)-x'(s)|^2ds,
\end{split}
\end{equation*}
and
\begin{equation*}
\begin{split}
&E\left|\varepsilon\int_0^tS_\varepsilon(t-s)g(x'(s),\mu'(s))d\omega(s)\right|^2\\
&\leqslant\varepsilon^2M^2e^{2\alpha T}(trQ)\int_0^tE|g(x'(s),\mu'(s))|^2ds\\
&\leqslant\varepsilon^2M^2e^{2\alpha T}(trQ)K_4\int_0^t(1+E|x'(s)|^2+\|\mu'(s)\|_\varphi^2)ds\\
&\leqslant\varepsilon^2M^2e^{2\alpha T}(trQ)K_4T[1+J(1+E|x_0|^2)+\mathop{{\rm sup}}_{t\in[0,T]}\|\mu'(t)\|_\varphi^2]\\
&<\infty.
\end{split}
\end{equation*}
By the Lebesgue's dominated convergence theorem, we obtain
$$E\left|\varepsilon\int_0^tS_\varepsilon(t-s)g(x'(s),\mu'(s))d\omega(s)\right|^2\rightarrow0\quad as\ \ \varepsilon\downarrow0$$
for any $\delta_3(\varepsilon)>0,$ there exists an $\varepsilon_3>\varepsilon\downarrow0$ and constant $l_3>0$ such that
$$E\left|\varepsilon\int_0^tS_\varepsilon(t-s)g(x'(s),\mu'(s))d\omega(s)\right|^2\leqslant l_3\delta_3(\varepsilon).$$
Let $\varepsilon_4=\min\{\varepsilon_1, \varepsilon_2, \varepsilon_3\}$, then for $\varepsilon_4>\varepsilon\downarrow0,$ we have
\begin{equation*}
\begin{split}
E|x_\varepsilon(t)-x'(t)|^2&\leqslant5\{l_1\delta_1(\varepsilon)+l_2\delta_2(\varepsilon)+l_3\delta_3(\varepsilon)\}\\
&\quad+20M^2e^{2\alpha T}[TK_1^2+\varepsilon^2(trQ)K_2^2]\int_0^tE|x_\varepsilon(s)-x'(s)|^2ds.
\end{split}
\end{equation*}
By the Gronwall's inequality, we have
$$E|x_\varepsilon(t)-x'(t)|^2\leqslant\tau(\varepsilon)\eta(t),\quad t\in[0,T],$$
where $\tau(\varepsilon)=5\{l_1\delta_1(\varepsilon)+l_2\delta_2(\varepsilon)+l_3\delta_3(\varepsilon)\}$;\ $\eta(t)=\exp\{20M^2e^{2\alpha T}[TK_1^2+\varepsilon^2(trQ)K_2^2]t\}$. The proof is completed.
\end{proof}
\section{Parametric Dependence Of The Equation}\label{sec5}

In this section, we will give a concrete application of Trotter-Kato approximation system introduced for the stochastic evolution equation, namely a classical limit theorem for the parameter dependence of System (\ref{system(1.1)})-(\ref{system(1.2)}).

First let us study the following family of stochastic evolution equations, which are the main objects of this section.
\begin{align}
\label{system(5.1)}
  dx_n(t)&=[A_nx_n(t)+f_n(x_n(t),\mu_n(t))]dt+g_n(x_n(t),\mu_n(t))d\omega(t),\quad t\in[0,T],\\
\label{system(5.2)}
  x_n(0)&=x_0,
 \end{align}
 where $A_n,\ n\in\mathbb{N^+},$ is the infinitesimal generator of a $C_0$-semigroup $\{S_n(t),\ t\in[0,T]\}$ of bounded linear operators on $H_1.$ $f_n,\ n\in\mathbb{N^+}$, is an appropriate $H_1$-valued function defined on $H_1\times (V_{\varphi^2}(H_1),\rho)$; $g_n,\ n\in\mathbb{N^+}$, is a $\mathcal{L}(H_2,H_1)$-valued function on $H_1\times (V_{\varphi^2}(H_1),\rho).$

 Let all terms in Equation (\ref{system(5.1)}) satisfy the conditions of Theorem (\ref{the3.1})(a). Then the equation has a unique mild solution $x_n$ in $C([0,T],L_2(\Omega;H_1)),$ for every $n\in\mathbb{N^+}$. Hence, we know $x_n(t)$ satisfies the following stochastic integral equation:
\begin{equation*}
\begin{split}
x_n(t)&=S_n(t)x_0+\int_0^tS_n(t-s)f_n(x_n(s),\mu_n(s))ds\\
&\quad+\int_0^tS_n(t-s)g_n(x_n(s),\mu_n(s))d\omega(s)\quad P-a.s.,
\end{split}
\end{equation*}
for $t\in[0,T].$ In order to study this equation and get the desired result, we need to make the following assumption.
\begin{spacing}{2.0}
\begin{large}
\noindent
\textbf{Assumption(A4)}
\end{large}
\end{spacing}

For each $n\in\mathbb{N^+}$, the nonlinear functions $f_n(x,\mu)$ and $g_n(x,\mu)$ satisfy the following Lipschitz and linear growth conditions for all $t\in[0,T]$:
$$|f_n(x,\mu)-f_n(y,\nu)| \leqslant K'_{1}(|x-y|+\rho(\mu,\nu)),\quad K'_1>0,\ x,y\in H_1,\ \mu,\nu\in V_{\varphi^2}(H_1),$$
$$|g_n(x,\mu)-g_n(y,\nu)| \leqslant K'_{2}(|x-y|+\rho(\mu,\nu)),\quad K'_2>0,\ x,y\in H_1,\ \mu,\nu\in V_{\varphi^2}(H_1),$$
$$|f_n(x,\mu)|^{p} \leqslant K'_{3}(1+|x|^{p}+\|\mu\|^{^{p}}_{\varphi}),\quad K'_3>0,\ x\in H_1,\ \mu\in V_{\varphi^2}(H_1),$$
$$|g_n(x,\mu)|^{p} \leqslant K'_{4}(1+|x|^{p}+\|\mu\|^{^{p}}_{\varphi}),\quad K'_4>0,\ x\in H_1,\ \mu\in V_{\varphi^2}(H_1),$$
for $p\geqslant2$. Note that the constants $K'_i,\ i=1,2,3,4$ don't depend on $t.$
\begin{spacing}{2.0}
\begin{large}
\noindent
\textbf{Assumption(A5)}
\end{large}
\end{spacing}

For every $Z>0$, we have the following two limits:
$$\mathop{{\rm lim}}_{n\rightarrow\infty}\mathop{{\rm sup}}_{|x|\leqslant Z}|f_n(x(t),\mu(t))-f(x(t),\mu(t))|\rightarrow0$$ and $$\mathop{{\rm lim}}_{n\rightarrow\infty}\mathop{{\rm sup}}_{|x|\leqslant Z}|g_n(x(t),\mu(t))-g(x(t),\mu(t))|\rightarrow0$$
uniformly in $t$ for all $t\in[0,T].$
\begin{theorem}
\label{the(5.1)}
Assume that (A1), (A2), (A4) and (A5) hold. Let $x(t)$ and $x_n(t)$ be the mild solutions of Equations (\ref{system(1.1)}) and (\ref{system(5.1)}), respectively. Then for each $T\in(0,\infty),$ we have
$$\mathop{{\rm sup}}_{0\leqslant t\leqslant T}E|x_n(t)-x(t)|^2\rightarrow0\quad n\rightarrow\infty.$$
\end{theorem}
\begin{proof}
It is clear that
\begin{equation*}
\begin{split}
&x_n(t)-x(t)\\
&=S_n(t)x_0-S(t)x_0+\int_0^t[S_n(t-s)f_n(x_n(s),\mu_n(s))-S(t-s)f(x(s),\mu(s))]ds\\
&\quad+\int_0^t[S_n(t-s)g_n(x_n(s),\mu_n(s))-S(t-s)g(x(s),\mu(s)]d\omega(s)\\
&=S_n(t)x_0-S(t)x_0+\int_0^t[S_n(t-s)[f_n(x_n(s),\mu_n(s))-f(x(s),\mu(s))]ds\\
&\quad+\int_0^t[S_n(t-s)[f_n(x_n(s),\mu_n(s))-f_n(x(s),\mu(s))]ds\\
&\quad+\int_0^t[S_n(t-s)-S(t-s)]f(x(s),\mu(s))ds\\
&\quad+\int_0^t[S_n(t-s)g_n(x_n(s),\mu_n(s))-S(t-s)g(x(s),\mu(s))]d\omega(s)\\
\end{split}
\end{equation*}
\begin{equation*}
\begin{split}
&\quad+\int_0^t[S_n(t-s)[g_n(x_n(s),\mu_n(s))-g_n(x(s),\mu(s))]d\omega(s)\\
&\quad+\int_0^t[S_n(t-s)-S(t-s)]g(x(s),\mu(s))d\omega(s),\quad P-a.s.
\end{split}
\end{equation*}
for all $t\in[0,T].$ Let the sum of the first term, the second term, the fourth term, the fifth term and the seventh term be abbreviated as $\beta(t).$ Now let us square this equation and then take the expectations on both sides, and we get
\begin{equation*}
\begin{split}
&E|x_n(t)-x(t)|^2\\
&\leqslant3\bigg\{E\left|\beta(t)\right|^2+E\left|\int_0^tS_n(t-s)[f_n(x_n(s),\mu_n(s))-f_n(x(s),\mu(s))]ds\right|^2\\
&\quad+E\left|\int_0^t[S_n(t-s)[g_n(x_n(s),\mu_n(s))-g_n(x(s),\mu(s))]d\omega(s)\right|^2\bigg\}\\
&\leqslant3\bigg\{E\left|\beta(t)\right|^2+TM^2e^{2\alpha T}\int_0^tE|f_n(x_n(s),\mu_n(s))-f_n(x(s),\mu(s))|^2ds\\
&\quad+(trQ)M^2e^{2\alpha T}\int_0^tE|g_n(x_n(s),\mu_n(s))-g_n(x(s),\mu(s))|^2ds\bigg\}\\
&\leqslant3\bigg\{E\left|\beta(t)\right|^2+TM^2e^{2\alpha T}{K'_{1}}^2\int_0^tE|x_n(s)-x(s)|^2+\rho^2(\mu_n(s),\mu(s))ds\\
&\quad+(trQ)M^2e^{2\alpha T}{K'_{2}}^2\int_0^tE|x_n(s)-x(s)|^2+\rho^2\mu_n(s),\mu(s))ds\bigg\}\\
&\leqslant3E\left|\beta(t)\right|^2+6M^2e^{2\alpha T}[T{K'_{1}}^2+(trQ){K'_{2}}^2]\int_0^tE|x_n(s)-x(s)|^2ds.
\end{split}
\end{equation*}
Then let $\theta=6M^2e^{2\alpha T}[T{K'_{1}}^2+(trQ){K'_{2}}^2],$ we obtain the following inequality from Lemma 1 in Gikhman and Skorohod \cite{GIISAV72}:
\begin{equation}
\begin{split}
\label{eq(5.3)}
E|x_n(t)-x(t)|^2&\leqslant3E|\beta(t)|^2+\theta\int_0^t\exp[\theta(t-s)]E|\beta(s)|^2ds.
\end{split}
\end{equation}
Hence, we only need to prove that $\mathop{{\rm sup}}_{0\leqslant t\leqslant T}E\left|\beta(t)\right|^2\rightarrow0$. The specific form of $E|\beta(t)|^2$ is as follows:
\begin{equation}
\begin{split}
\label{eq(5.4)}
E|\beta(t)|^2&\leqslant5\bigg\{E|S_n(t)x_0-S(t)x_0|^2+E\left|\int_0^tS_n(t-s)[f_n(x(s),\mu(s))-f(x(s),\mu(s))]ds\right|^2\\
&\quad+E\left|\int_0^t[S_n(t-s)-S(t-s)]f(x(s),\mu(s))ds\right|^2\\
&\quad+E\left|\int_0^tS_n(t-s)[g_n(x(s),\mu(s))-g(x(s),\mu(s))]d\omega(s)\right|^2\\
&\quad+E\left|\int_0^t[S_n(t-s)-S(t-s)]g(x(s),\mu(s))d\omega(s)\right|^2\bigg\}.
\end{split}
\end{equation}
Next, we show that each term of Equation (\ref{eq(5.4)}) tends to zero as $n\rightarrow\infty.$ The first term has been proved in Equation (\ref{eq(4.8)}). Thus we first study the second term, by Assumption (A1), (A4) and Theorem \ref{the3.1}, we have
\begin{equation*}
\begin{split}
&E\left|\int_0^tS_n(t-s)[f_n(x(s),\mu(s))-f(x(s),\mu(s))]ds\right|^2\\
&\leqslant TM^2e^{2\alpha T}\int_0^tE\left|f_n(x(s),\mu(s))-f(x(s),\mu(s))\right|^2ds\\
&\leqslant2TM^2e^{2\alpha T}\int_0^tE[|f_n(x(s),\mu(s))|^2+|f(x(s),\mu(s))|^2]ds\\
&\leqslant2TM^2e^{2\alpha T}(K'_{3}+K_3)\int_0^t(1+E|x(s)|^2+\|\mu(s)\|_\varphi^2)ds\\
&\leqslant2T^2M^2e^{2\alpha T}(K'_{3}+K_3)[1+J(1+E|x_0|^2)+\mathop{{\rm sup}}_{0\leqslant t\leqslant T}\|\mu(t)\|_\varphi^2]\\
&<\infty.
\end{split}
\end{equation*}
Therefore, by Assumption (A5) and the Lebesgue's dominated convergence theorem, we obtain
$$\mathop{{\rm lim}}_{n\rightarrow\infty}\mathop{{\rm sup}}_{0\leqslant t\leqslant T}E\left|\int_0^tS_n(t-s)[f_n(x(s),\mu(s))-f(x(s),\mu(s))]ds\right|^2=0.$$
From the Assumption (A1) and (A2), we get
\begin{equation*}
\begin{split}
&E\left|\int_0^t[S_n(t-s)-S(t-s)]f(x(s),\mu(s))ds\right|^2\\
&\leqslant T^2M^2e^{2\alpha T}K_3[1+J(1+E|x_0|^2)]\mathop{{\rm sup}}_{0\leqslant t\leqslant T}\|\mu(t)\|_\varphi^2]\\
&<\infty.
\end{split}
\end{equation*}
Hence, by Theorem \ref{the4.1} (b) and the Lebesgue's dominated convergence theorem, we have
$$\mathop{{\rm lim}}_{n\rightarrow\infty}\mathop{{\rm sup}}_{0\leqslant t\leqslant T}E\left|\int_0^t[S_n(t-s)-S(t-s)]f(x(s),\mu(s))ds\right|^2=0.$$
Then we consider the remaining two stochastic integral terms, by Assumption (A1), (A4) and the Proposition \ref{por(2.3)}, we have
\begin{equation*}
\begin{split}
&E\left|\int_0^tS_n(t-s)[g_n(x(s),\mu(s))-g(x(s),\mu(s))]d\omega(s)\right|^2\\
&\leqslant (trQ)M^2e^{2\alpha T}\int_0^tE\left|g_n(x(s),\mu(s))-g(x(s),\mu(s))\right|^2ds\\
&\leqslant2(trQ)M^2e^{2\alpha T}\int_0^tE[|g_n(x(s),\mu(s))|^2+|g(x(s),\mu(s))|^2]ds\\
&\leqslant2(trQ)M^2e^{2\alpha T}(K'_{4}+K_4)\int_0^t(1+E|x(s)|^2+\|\mu(s)\|_\varphi^2)ds\\
&\leqslant2(trQ)M^2e^{2\alpha T}(K'_{4}+K_4)T[1+J(1+E|x_0|^2)+\mathop{{\rm sup}}_{0\leqslant t\leqslant T}\|\mu(t)\|_\varphi^2]\\
&<\infty.
\end{split}
\end{equation*}
Thus, by Assumption (A5) and the Lebesgue's dominated convergence theorem, we have
$$\mathop{{\rm lim}}_{n\rightarrow\infty}\mathop{{\rm sup}}_{0\leqslant t\leqslant T}E\left|\int_0^tS_n(t-s)[g_n(x(s),\mu(s))-g(x(s),\mu(s))]d\omega(s)\right|^2=0.$$
Next we consider the second stochastic term,
\begin{equation*}
\begin{split}
&E\left|\int_0^t[S_n(t-s)-S(t-s)]g(x(s),\mu(s))d\omega(s)\right|^2\\
&\leqslant (trQ)M^2e^{2\alpha T}K_4[1+J(1+E|x_0|^2)]\mathop{{\rm sup}}_{0\leqslant t\leqslant T}\|\mu(t)\|_\varphi^2]\\
&<\infty.
\end{split}
\end{equation*}
Then by the Theorem \ref{the4.1} (b) and the Lebesgue's dominated convergence theorem, we have
$$\mathop{{\rm lim}}_{n\rightarrow\infty}\mathop{{\rm sup}}_{0\leqslant t\leqslant T}E\left|\int_0^t[S_n(t-s)-S(t-s)]g(x(s),\mu(s))d\omega(s)\right|^2=0.$$
Finally, we get
$$\mathop{{\rm lim}}_{n\rightarrow\infty}\mathop{{\rm sup}}_{0\leqslant t\leqslant T}E|\beta(t)|^2=0.$$
Substituting the above equation into Equation (\ref{eq(5.3)}), we obtain
$$\mathop{{\rm lim}}_{n\rightarrow\infty}\mathop{{\rm sup}}_{0\leqslant t\leqslant T}E|x_n(t)-x(t)|^2=0.$$
The proof is completed.
\end{proof}
\begin{corollary}
Supposed that the coefficients in the Equation (\ref{system(5.1)}) depend on the parameter $\lambda$ varying by a set of numbers $N$.
\begin{align}
\label{system(5.5)}
  dx_\lambda(t)&=[A_\lambda x_\lambda(t)+f_\lambda(x_n(t),\mu_\lambda(t))]dt+g_\lambda(x_\lambda(t),\mu_\lambda(t))d\omega(t),\quad t\in[0,T],\\
\label{system(5.6)}
  x_\lambda(0)&=x_0,
 \end{align}
 where $A_\lambda$ is the infinitesimal generator of a $C_0$-semigroup $\{S_\lambda(t),t\geqslant0\}$ of bounded linear operators on $H_1.$ $f_\lambda$ is an appropriate $H_1$-valued function defined on $H_1\times (V_{\varphi^2}(H_1),\rho);$ $g_\lambda$ is a $\mathcal{L}(H_2,H_1)$-valued function on $H_1\times (V_{\varphi^2}(H_1),\rho).$ Then we need to make similar assumption.
\begin{spacing}{2.0}
\begin{large}
\noindent
\rm{\textbf{Assumption(A6)}}
\end{large}
\end{spacing}

For each $Z>0,$ we have the following two limits:
$$\mathop{{\rm lim}}_{\lambda\rightarrow\lambda'}\mathop{{\rm sup}}_{|x|\leqslant Z}|f_\lambda(x(t),\mu(t))-f_\lambda'(x(t),\mu(t))|\rightarrow0$$ and $$\mathop{{\rm lim}}_{\lambda\rightarrow\lambda'}\mathop{{\rm sup}}_{|x|\leqslant Z}|g_\lambda(x(t),\mu(t))-g_\lambda'(x(t),\mu(t))|\rightarrow0$$
uniformly in $t$ for each $t\in[0,T].$

Let $A_\lambda,\ A_{\lambda'}$ be infinitesimal generators of $C_0$-semigroups $S_\lambda(t)$ and $S_\lambda'(t)$, respectively, and $D(A_\lambda)=D(A_{\lambda'}),$ where $\lambda\in N$ and $\mathop{{\rm lim}}_{\lambda\rightarrow\lambda'}S_\lambda(t)=S_{\lambda'}(t)$ uniformly in $t\in[0,T]$ for every $T\in(0,\infty).$ Finally, let $A_\lambda,$ $f_\lambda(x,\mu)$ and $g_\lambda(x,\mu)$ satisfy the conditions of Theorem \ref{the3.1} and appropriate assumption for each $\lambda\in N.$ Then the Equation (\ref{system(5.5)}) has a unique mild solution $x_\lambda(t)$ and satisfies
$$\mathop{{\rm lim}}_{\lambda\rightarrow\lambda'}\mathop{{\rm sup}}_{t\in[0,T]}E|x_\lambda(t)-x_{\lambda'}(t)|^2=0.$$
\end{corollary}
\begin{proof}
The proof can be obtained directly from the proof of Theorem \ref{the(5.1)}.
\end{proof}

\section{An Example}

This example is adapted from Kunze and van Neerven \cite{KMVN11}. Let us first study the stochastic partial differential equation as follows:
\begin{equation}
\begin{split}
\label{system(6.1)}
\frac{\partial x}{\partial t}x(t,z)&=[Ax(t,z)+f(x(t,z),\mu(t))]dt+\sum_{k=1}^{K}g_k(x(t,z),\mu(t))\frac{\omega_k}{\partial t}(t),\quad z\in O,\ t\in[0,T],\\
x(t,z)&=0,\quad z\in\partial O,\ t\in[0,T],\\
x(0,z)&=x_0,\quad z\in O,
\end{split}
\end{equation}
where $\{\omega_k(t),\ k\in\mathbb{N^+}\}$ are independent real-valued standard Wiener processes. $O$ is a bounded open field in $R^m$. $A$ is the second-order divergence form differential operator defined by
$$Ax=\sum_{i=1}^m\frac{\partial}{\partial z_i}\left(q_{ij}(z)\sum_{j=1}^{m}\frac{\partial x}{\partial z_j}(z)\right)+\sum_{j=1}^mr_j(z)\frac{\partial x}{\partial z_j}(z),$$
let the parameters $q=(q_{ij})$ and $r=(r_j)$ satisfy suitable boundedness and uniform ellipticity conditions. In Equation (\ref{system(6.1)}), the functions $f$ and $g_k$ are Lipschitz continuous. Let the following expression denote the Lipschitz seminorm of function $f$,
$$\|f\|_{Lip}=\mathop{{\rm sup}}_{u\neq v}\frac{|f(u)-f(v)|}{|u-v|}.$$

In order to get the desired results, we need to make the following assumptions.
\begin{spacing}{2.0}
\begin{large}
\noindent
\textbf{Assumption(A7)}
\end{large}
\end{spacing}

Let $q,$ $q_n\in L_\infty(O;R^{m\times m})$ and $r,$ $r_n\in L_\infty(O;R^m).$ Let $f,\ f_n,\ g_k,\ {(g_k)}_n$ be Lipschitz continuous. Suppose that there exist finite numbers $\sigma$ and $N$ such that
\begin{spacing}{1.5}
\end{spacing}
\noindent
(i) $q,$ $q_n$ are symmetric and $qx\cdot x,$ $q_nx\cdot x\geqslant\sigma|x|^2$ for all $x\in R^m;$

\noindent
(ii) $\|q\|_\infty,$ $\|q_n\|_\infty,$ $\|r\|_\infty,$ $\|r_n\|_\infty\leqslant N;$

\noindent
(iii) $\|f\|_{Lip},$ $\|f_n\|_{Lip},$ $\|g_k\|_{Lip},$ $\|{(g_k)}_n\|_{Lip}\leqslant N;$

\noindent
(iv) $q_n=q,$ $r_n=r$ as $n\rightarrow\infty$ a.e. on $O;$

\noindent
(v) $f_n=f,$ ${(g_k)}_n=g_k$ as $n\rightarrow\infty$ pointwise on $O.$
\begin{spacing}{2.0}
\begin{large}
\noindent
\textbf{Assumption(A8)}
\end{large}
\end{spacing}

(i) The operators $A$ and $A_n,$ for each $n\in\mathbb{N^+}$ are densely defined, closed and uniformly
\par\setlength\parindent{1.5em}
sectorial on $H_1$ in the sense, there exists constants $M\geqslant1$ and $\alpha\in R$ such that $A$ and
\par\setlength\parindent{1.5em}
$A_n$ is sectorial of type $(M,\alpha).$

(ii) The operators $A_n$ converge to $A$ in the strong resolvent sense:
$$\mathop{{\rm lim}}_{n\rightarrow\infty}R(\lambda,A_n)x=R(\lambda,A)x$$

for all $Re\lambda>\alpha$ and $x\in H_1.$

Under the Assumption (A8), the operators $A_n$ and $A$ generate strongly continuous analytic semigroups $S(t)$ and $S_n(t)$, respectively, satisfying the following uniform bounds:
$$\|S(t)\|,\ \|S_n(t)\|\leqslant Me^{\alpha t},\quad t\in[0,T],$$
$$\|AS(t)\|,\ \|A_nS_n(t)\|\leqslant\frac{M'}{t}e^{\alpha t},\quad t\in[0,T].$$

In order to reconstruct the above equation under the abstract background of Banach space $L_p(O),$ $1<p<\infty,$ we use the variational approach. Consider the sesquilinear form as follows:
$$q[u,v]:=\int_O(p\nabla u)\cdot\overline{\nabla v}+(r\cdot\nabla u)\overline{v}dx$$
on the domain $D(q):=H_0^1(O).$ The sectorial operator $A$ on $L_p(O)$ associated with $q$ generates a strongly continuous analytic semigroup $\{S(t),\ t\in[0,T]\}$ and extrapolates it to a consistent family of strongly continuous analytic semigroups $\{S^{(p)}(t),\ t\in[0,T])\}$ on $L_p(O)$. We use $A^{(p)}$ to represent its corresponding infinitesimal generator. Similarly, the form $q_n$ and the semigroups $S_n^{(p)}(t)$ related to the generators $A_n^{(p)}$ are defined.
\begin{proposition}
\label{por(6.1)}
\rm(Kunze and van Neerven \cite{KMVN11}) \it Let the Assumption (A7) (i), (ii) and (iv) hold. Then the operators $A^{(p)}$ and $A_n^{(p)}$ satisfy the conditions of Assumption (A8).
\end{proposition}
\begin{proposition}
\label{por(6.2)}
\rm(Kunze and van Neerven \cite{KMVN11}) \it Assume that Assumption (A7) (iii) and (v) hold. Then

\noindent
(i) the maps $f,\ f_n:\ L_p(O)\times (V_{\varphi^2}(L_p(O)),\rho)\rightarrow L_p(O)$ defined by
$$[f(x,\mu)](z):=f(x(z),\mu),\ \ [f_n(x,\mu)](z):=f_n(x(z),\mu),$$
where $\mu\in C([0,T],(V_{\varphi^2}(L_p(O)),\rho))$, satisfy the conditions of Assumption (A1) and (A4),

\noindent
(ii) the maps $g,\ g_n:\ L_p(O)\times (V_{\varphi^2}(L_p(O)),\rho)\rightarrow\mathcal{L}(R_K, L_p(O))$ defined by
$$[g(x,\mu)h](z):=\sum_{k=1}^Kg_k(x(z),\mu)(e_k,h),\ \ [g_n(x,\mu)h](z):=\sum_{k=1}^K{(g_k)}_n(x(z),\mu)(e_k,h),$$
where $\{e_k\}_{k=1}^K$ is the standard unit basis of $R^K$ and $\mu\in C([0,T],(V_{\varphi^2}(L_p(O)),\rho))$, satisfy the conditions of Assumption (A1) and (A4).
\end{proposition}

Hence, System (\ref{system(6.1)}) can be written in the abstract form as System (\ref{system(1.1)})-(\ref{system(1.2)}). Therefore, it is easy to verify that the classical limit theorem, as a concrete application, is also true for the stochastic partial differential equation in the form of System (\ref{system(6.1)}).

\bigskip

\noindent{\bf Acknowledgments:}

The study was supported by the Humanities and Social Science Foundation of Ministry of Education (Grant No. 20YJC790174), the Natural Science Foundation of Tianjin city (Grant No. 18JCYBJC18900) and the National Natural Science Foundation of China (Grant No. 11301380).


\end{document}